\pdfoutput=1
\RequirePackage{fix-cm}
\documentclass[reqno]{amsproc}
\usepackage{amssymb}
\usepackage{hyperref}
\usepackage{siunitx}
\usepackage{longtable}

\usepackage{tikz}

\usepackage{mathtools}
\mathtoolsset{showonlyrefs}

\usepackage{listings}
\usepackage[shortalphabetic,nobysame]{amsrefs}
\usepackage{xyzbib}

\numberwithin{equation}{section}

\newtheorem{theorem}{Theorem}
\newtheorem{proposition}[theorem]{Proposition}
\newtheorem{conjecture}[theorem]{Conjecture}

\newtheorem*{remark}{Remark}

\let\leq=\leqslant
\let\geq=\geqslant

\newcommand{\rmi}{\mathrm{i}}
\newcommand{\rmd}{\mathrm{d}}

\newcommand{\yc}{y_\mathrm{c}}

\newcommand{\rowspace}{\\[3pt]}

\begin{document}
\title{Frozen-corner enumeration of Alternating Sign Matrices}

\author{Filippo Colomo}
\address{INFN, Sezione di Firenze\\
Via G. Sansone 1, 50019 Sesto Fiorentino (FI), Italy}
\email{colomo@fi.infn.it}

\author{Andrei G. Pronko}
\address{Steklov Mathematical Institute, 
Fontanka 27, 191023 Saint Petersburg, Russia}
\email{a.g.pronko@gmail.com}

\begin{abstract}
An Alternating Sign Matrix (ASM) is a square matrix with entries in
$\{0,1,-1\}$, and such that: $i)$ in each row and columns, nonzero
entries alternate in sign; $ii)$ for any given row or column, entries
sum up to $1$.  We define the frozen-square enumeration as the
enumeration of $n\times n$ ASMs under the refinement of having,
located in a corner, an $s\times s$ square of entries that are all
zeroes.  We state a conjectural formula for such enumeration, in terms of
the determinant of some $s\times s$ matrix whose entries are given explicitly. 
We provide numerical support in favour of our conjecture. We also
illustrate the relevance of the conjectured formula in connection with
the limit shape observed in large ASMs, its fluctuations, and the
Tracy--Widom distribution.
\end{abstract}

\maketitle
\section{Introduction}

\subsection{Alternating Sign Matrices}
An Alternating Sign Matrix (ASM) is a square matrix with entries in
$\{0,1,-1\}$, and such that: $i)$ in each row and columns, nonzero
entries alternate in sign; $ii)$ for any given row or column, entries
sum up to $1$.  ASMs were introduced by Mills, Robbins, and Rumsey in
connections with descending plane partitions
\cite{MRR-82,MRR-83}. There are many interesting results concerning
ASMs: for an introductory review, see the nice book by Bressoud
\cite{Br-99}, while, for recent developments, see, e.g., \cite{BFK-23} and
reference therein.

Let $A_n$ denote the number of $n\times n$ ASMs. It was conjectured in
\cite{MRR-83}, and proved by Zeilberger \cite{Ze-94}, that
\begin{equation}\label{eq:asm}
A_n=\prod_{j=0}^{n-1}\frac{(3j+1)!}{(n+j)!}.
\end{equation}
An alternative proof was provided soon after by Kuperberg
\cite{K-96}. It is based on  a simple bijection, first described
by Robbins and Rumsey \cite{RR-86}, see also \cite{EKLP-92}, between
ASMs and configurations of square ice \cite{Li-67a}. The proof uses
Izergin--Korepin formula for the square ice partition function
\cite{K-82,I-87}, whose derivation relies on the integrability of the
model \cite{B-82} and the Quantum Inverse Scattering Method
\cite{TF-79,KBI-93}.

Let $A_{n,r}$, $r=1,\dots,n$, denote the number of $n\times n$ ASMs
whose sole nonzero entry, equal to $1$, in the first (or last) row (or
column) occurs at the $r$th position. This is referred to as the
refined enumeration of ASMs.  It was conjectured in \cite{MRR-83}, and
proven by Zeilberger \cite{Ze-96}, that
\begin{equation}\label{eq:refined}
  A_{n,r}=\binom{n+r-2}{r-1}\frac{(2n-r-1)!}{(n-r)!}
  \prod_{j=0}^{n-2}\frac{(3j+1)!}{(n+j)!}.
\end{equation}
The proof is built on Kuperberg's idea of using the bijection with
square ice and Izergin--Korepin formula. Alternative proofs of
\eqref{eq:asm} and \eqref{eq:refined} have been proposed over the
years, see e.g. \cite{CP-05b,F-07,F-16,FK-20,FK-22}.

ASMs, and square ice, can conveniently be described in terms of paths
on a square grid.  More precisely, let us consider $n$ paths on the
$n\times n$ square grid. Allowed path configurations are such that
each vertex of the grid may be in one out of six allowed
configurations, that we call of type $1, \dots,6$, respectively, see
Fig.~\ref{fig-vertices_entries}.  It follows that the vertices of type
5 and 6 (turning vertices) must alternate along each horizontal or
vertical line of the grid. We require the $n$ paths to be densely
packed on the top and left boundary of the grid, see
Fig.~\ref{fig-dwbc_asm_mt}; this is known as the `Domain Wall'
boundary condition. The imposed boundary condition implies
that the first and last `turning' vertices on each horizontal or
vertical line must all be of type 5.  The bijection with $n\times n$
ASMs follows.

\begin{figure}
\centering

\begin{tikzpicture}[scale=.5]
\draw  (0,1)--(2,1);
\draw  (1,0)--(1,2);
\draw  [very thick] (3,1)--(5,1);
\draw  [very thick] (4,0)--(4,2);
\draw  (6,1)--(8,1);
\draw  [very thick] (7,0)--(7,2);
\draw  (10,0)--(10,2);
\draw  [very thick] (9,1)--(11,1);
\draw  (13,0)--(13,1)--(14,1);
\draw  [very thick] (12,1)--(13,1)--(13,2);
\draw  (15,1)--(16,1)--(16,2);
\draw  [very thick] (16,0)--(16,1)--(17,1);
\node  at (1,-1) {$0$};
\node  at (4,-1) {$0$};
\node  at (7,-1) {$0$};
\node  at (10,-1) {$0$};
\node  at (13,-1) {$1$};
\node  at (16,-1) {$-1$};
\end{tikzpicture}

\caption{The six allowed vertex configurations
for the square ice in the path description
and the corresponding entries for ASMs}
\label{fig-vertices_entries}
\end{figure}

\begin{figure}
\centering 

\begin{tikzpicture}[scale=.5]
\draw (1,0) -- (1,6);
\draw (2,0) -- (2,6);
\draw (3,0) -- (3,6);
\draw (4,0) -- (4,6);
\draw (5,0) -- (5,6);
\draw (0,1) -- (6,1);
\draw (0,2) -- (6,2);
\draw (0,3) -- (6,3);
\draw (0,4) -- (6,4);
\draw (0,5) -- (6,5);
\draw [ultra thick] (0,5)--(1,5)--(1,6);
\draw [ultra thick] (0,4)--(1,4)--(1,5)--(2,5)--(2,6);
\draw [ultra thick] (0,3)--(1,3)--(1,4)--(2,4)--(2,5)--(3,5)--(3,6);
\draw [ultra thick] (0,2)--(2,2)--(2,4)--(3,4)--(3,5)--(4,5)--(4,6);
\draw [ultra thick] (0,1)--(3,1)--(3,2)--(4,2)--(4,4)--(5,4)--(5,6);

\node at (11,3)
{\setlength{\arraycolsep}{5pt}
$\begin{pmatrix} 
0 & 0 & 0 & 1 & 0 \\[3pt] 
0 & 0 & 1 & -1 & 1 \\[3pt] 
1 & 0 & 0 & 0 & 0 \\[3pt] 
0 & 1 & -1 & 1 & 0 \\[3pt] 
0 & 0 & 1 & 0 & 0 
\end{pmatrix}$};

\node at (20,3) 
{$\begin{matrix} 4 \\[3pt]
3 \qquad 5 \\[3pt] 
1 \qquad 3 \qquad  5\\[3pt]
1 \qquad 2 \qquad 4 \qquad 5\\[3pt]
1 \qquad 2 \qquad 3 \qquad 4 \qquad 5
\end{matrix}$};
\end{tikzpicture}

\caption{A path configuration on the $5\times 5$ grid with Domain Wall boundary
conditions (left), the corresponding ASM (center), and the related
monotone triangle (right).}
\label{fig-dwbc_asm_mt}
\end{figure}

Under this bijection, the refined enumeration of $n\times n$ ASMs
becomes the enumeration of path configurations on the $n\times n$
square grid, conditioned to have their sole vertex of type 5 at the
$r$th position in the first (or last) horizontal (or vertical)
line. Focussing of the last column for definiteness, see
Fig.~\ref{fig-refined}, it is clear that the refined enumeration is
also equivalent to the enumeration of path configurations on the
$n\times (n-1)$ square grid, with $n-1$ paths flowing vertically from
the top boundary, and one horizontally, from the right boundary, at
the $r$th position.

\begin{figure}
\centering

\begin{tikzpicture}[scale=.5]
\node at (-5,3) {
$\begin{pmatrix} 
* & * & * & * & 0 \\[3pt]
* & * & * & * & 0 \\[3pt]
* & * & * & * & 0 \\[3pt]
* & * & * & * & 1 \\[3pt]
* & * & * & * & 0 \\[3pt]
\end{pmatrix}
$};

\draw [help lines] (0,1)--(6,1);
\draw [help lines] (0,2)--(6,2);
\draw [help lines] (0,3)--(6,3);
\draw [help lines] (0,4)--(6,4);
\draw [help lines] (0,5)--(6,5);
\draw [help lines] (1,0)--(1,5);
\draw [help lines] (2,0)--(2,5);
\draw [help lines] (3,0)--(3,5);
\draw [help lines] (4,0)--(4,5);
\draw [help lines] (5,0)--(5,5);
\draw [very thick](1,5)--(1,6);
\draw [very thick](2,5)--(2,6);
\draw [very thick](3,5)--(3,6);
\draw [very thick](4,5)--(4,6);
\draw [very thick](5,5)--(5,6);
\draw [very thick](0,1)--(1,1);
\draw [very thick](0,2)--(1,2);
\draw [very thick](0,3)--(1,3);
\draw [very thick](0,4)--(1,4);
\draw [very thick](0,5)--(1,5);
\draw [very thick] (4,2)--(5,2)--(5,5);
\fill [nearly transparent,gray] (4.2,.8)--(4.2,5.2)--(.8,5.2)--(.8,.8);
\draw [|->] (6.5,5.3)--(6.5,1.7);
\node at (7,3.5) {$r$};
\draw [help lines] (9,1)--(14,1);
\draw [help lines] (9,2)--(14,2);
\draw [help lines] (9,3)--(14,3);
\draw [help lines] (9,4)--(14,4);
\draw [help lines] (9,5)--(14,5);
\draw [help lines] (10,0)--(10,5);
\draw [help lines] (11,0)--(11,5);
\draw [help lines] (12,0)--(12,5);
\draw [help lines] (13,0)--(13,5);
\draw [very thick](10,5)--(10,6);
\draw [very thick](11,5)--(11,6);
\draw [very thick](12,5)--(12,6);
\draw [very thick](13,5)--(13,6);
\draw [very thick](9,1)--(10,1);
\draw [very thick](9,2)--(10,2);
\draw [very thick](9,3)--(10,3);
\draw [very thick](9,4)--(10,4);
\draw [very thick](9,5)--(10,5);
\draw [very thick] (13,2)--(14,2);
\fill [nearly transparent,gray] (13.2,.8)--(13.2,5.2)--(9.8,5.2)--(9.8,.8);
\end{tikzpicture}

\caption{Refined enumeration of ASMs: 
A generic ASM with the sole nonzero entry of the last column at
position $r$ (left), the corresponding path configuration on the
square grid (center), and the modified grid with its last line removed
(right). Path configurations within the shaded area are generic.}
\label{fig-refined}
\end{figure}

\subsection{Frozen square enumeration of ASMs}
In the present paper we are interested in the enumeration of ASMs
under a different refinement, namely that of having in, say, its
top-left corner, a square of entries that are all zeroes (a frozen
square). More precisely, a given ASM $(a_{ij})_{1\leq i,j\leq n}$ has
an $s\times s$ frozen square in its top-left corner if $a_{ij}=0$ for
all $i,j=1,\dots,s$.  Let $B_{n,s}$, $s=1,\dots,n$, be the number of
such matrices.  We call the quantity $B_{n,s}$ \emph{frozen-corner}
enumeration. Clearly, by symmetry, $B_{n,s}$ does not depend on which
of the four corners we choose to freeze.

Let us briefly mention some simple properties of the quantity
$B_{n,s}$.  When $s=1$, a zero is required to be in the top left
corner. Hence, the sole $1$ of the first row must be at position
$r=2,\dots,n$. It follows that
\begin{equation}\label{eq:property1}
B_{n,1}=\sum_{r=2}^n A_{n,r}=\sum_{r=1}^{n-1} A_{n,r},
  \end{equation}
where the second equality follows from the set of identities
$A_{n,r}=A_{n,n-r+1}$.  Further properties of the frozen-corner
enumeration may be easily derived by using  the path description
of square ice.  For instance, we have
\begin{equation}\label{eq:property2}
  B_{n,s}=0, \qquad  s>\lfloor n/2\rfloor.
\end{equation}
We also have, for ASMs of even size:
\begin{equation}\label{eq:property3}
  B_{2s,s}=(A_s)^2,
\end{equation}
and for ASMs of odd size:
\begin{equation}\label{eq:property4}
  B_{2s+1,s}=\sum_{j,k=1}^{s+1}
\binom{j+k-2}{j-1}  A_{s+1,j}A_{s+1,k}.
\end{equation}
The derivation of this last identity is slightly more involved. We
sketch it in Appendix \ref{app:B_proof}.

\subsection{A conjecture}
Before stating our conjecture for the numbers $B_{n,s}$, it is
convenient to introduce some notations.
Let us introduce the (normalized) generating function for the refined
enumeration of ASMs \eqref{eq:refined}:
\begin{equation}\label{eq:gnz}
h_n(z):=\frac{1}{A_n}\sum_{r=1}^n A_{n,r}z^{r-1}.
\end{equation}
For $i= 1,2,\dots$, let
\begin{equation}\label{eq:fpm}
  f_i^{\pm}(z):=\left[1\pm(-1)^i z\right] \frac{(1-z)^{i-1}}{z^i}
        h_{n-s+i}(z).
\end{equation}
For any given $n$ and $s$, let us introduce the $s\times s $ matrix
        $M$, with entries
\begin{equation}\label{eq:matrix}
M_{ij}:=
\frac{A_{n-s+j}}{A_{n-s+j-1}}\oint_{C_0}\oint_{C_0}\frac{f_i^+(z)f_j^-(w)}{1-z-w}
\frac{\rmd z\rmd w}{(2\pi\rmi)^2} ,\qquad i,j=1,\dots, s,
\end{equation}
where $C_0$ is a small closed contour around the origin, containing
no other singularity of the integrand.

\begin{remark}
  The entries of the $s\times s$ matrix $M$ may be expressed in terms
  of ASM's plain and refined enumerations, as follows:
  \begin{multline}\label{eq:matrix_new}
    M_{ij}=\frac{1}{A_{n-s+i}A_{n-s+j-1}}
    \sum_{k\geq 0}\sum_{l=0}^k 
    \sum_{p=1}^{n-s+i} \sum_{q=1}^{n-s+j}  (-1)^{i+j-k-p-q}\binom{k}{l} 
    \\
    \times \left[\binom{i-1}{i+l-k-p} -(-1)^i \binom{i-1}{i+l-k-p-1}\right]
    \\ \times
    \left[\binom{j-1}{j-l-q} +(-1)^j\binom{j-1}{j-l-q-1}\right]
    A_{n-s+i,p}A_{n-s+j,q} 
\end{multline} with $i,j=1,\dots,s$. This
  is easily worked out by
  Laurent expanding  
  the integrand in \eqref{eq:matrix} in the
  variables $z$ and $w$, and evaluating the corresponding residues at
  the origin (that is, the  constant term).
\end{remark}

We propose the following formula for the frozen-corner enumeration.
\begin{conjecture}\label{main_conj} The number $B_{n,s}$ may be evaluated
  in terms of a
  determinant:
  \begin{equation}
    B_{n,s}=A_n \det_{1\leq i,j\leq s}(1-M)
  \end{equation}
  where the $s\times s$ matrix $M$ is defined in
  \eqref{eq:matrix_new}, or equivalently, in terms of
  Eqs. \eqref{eq:gnz}, \eqref{eq:fpm}, and \eqref{eq:matrix}.
\end{conjecture}

Conjecture \ref{main_conj} has already been stated in the physical
context of the six-vertex model \cite{CP-24}.  
Here we reformulate the assertion in the ASM language,
to bring it to the attention of the combinatorics community, where
somebody might be able to build a rigorous proof. We also provide
strong numerical evidence in support of it.

In Section \ref{sec:derivation} we explain the rationale behind
Conjecture \ref{main_conj}.

In Section \ref{sec:numerics} we show how the numbers $B_{n,s}$ may be
evaluated, for moderate values of $n$ ($n\leq 20$), by direct
enumeration of corresponding monotone triangles satisfying some simple
constraint. The obtained values coincides with those provided by
evaluating the right hand side of Conjecture \ref{main_conj}.

In Section \ref{sec:context} we provide further motivations and
context behind the work that led us to state Conjecture
\ref{main_conj}. We also illustrate some nontrivial implications of it,
in connection with the limit shape observed in large ASMs, its
fluctuations, and the Tracy--Widom distribution.

\section{Rationale}\label{sec:derivation}

\subsection{Why the frozen-square enumeration?}
The motivation in studying frozen-corner enumeration of ASMs stems
from asymptotic combinatorics. Large ASMs are known to develop a limit
shape, with four `frozen' regions, containing only zeroes, in the four
corners, and a central `disordered' region, containing all three kind
of entries. The phenomenon is closely related to the Arctic Circle
Theorem \cite{JPS-98}. Indeed this occurs in the domino tilings of the
Aztec Diamond, which can be mapped onto the 2-enumeration of ASMs
\cite{EKLP-92}, i.e., when a weight $2^k$ is given to each ASM with
$k$ `$-1$' entries \cite{MRR-82}.

In order to study the phenomenon, one could try to build a closed form
expression for ASM enumeration under the refinement that, say,
$a_{ij}=-1$ for some given $i,j$. Such an enumeration is currently
beyond our capabilities.  A simpler refinement consists in requiring
the presence of a frozen square (or, more generally, rectangle) in one
of the corners. This is the frozen-square enumeration described above.
As we shall see, this may yield relevant information on the
limit-shape phenomenon in large ASMs.

\subsection{A closed formula for  the frozen-square enumeration}
Over the years, various closed form formulae have been provided for
the frozen-square enumeration \cite{CP-07b,CDP-21}. Here we present
the most relevant one for our purposes.
\begin{proposition}\label{MIR}
The frozen-square enumeration may be evaluated in terms of the
multiple integral
\begin{multline}\label{eq:efpMIR}
  B_{n,s}= (-1)^s A_n \oint_{C_0}\dots \oint_{C_0}
  \prod_{j=1}^s\frac{1}{z_j^{n-s}(z_j-1)^{s-j+1}}
  \\ \times\prod_{1\leq j<k\leq s}
  \frac{1}{z_jz_k-z_j+1}\det_{1\leq j,k\leq s}
  \left[(z_j-1)^{s-k}z_j^{k-1}h_{n-k+1}(z_j)\right]
  \frac{\rmd^sz}{(2\pi\rmi)^s},
\end{multline}
where $h_n(z)$ is defined in \eqref{eq:gnz}, and $C_0$ is a small
closed contour around the origin.
\end{proposition}
Relation \eqref{eq:efpMIR} is essentially a `constant term identity'
for $B_{n,s}$. Proposition \ref{MIR} is a particular case of a more
general result, valid when the frozen region is a rectangle, and
originally derived in the context of the six-vertex model, for
arbitrary values of its parameters, using integrability
techniques \cite{B-82}, and the Quantum Inverse Scattering
Method \cite{TF-79,KBI-93}. Relation \eqref{eq:efpMIR} may be obtained
from, e.g., Eq.~(5.14) in \cite{CDP-21}, simply by replacing therein
$N\to n$, $r\to n-s$, specializing the values of the parameters to
$\Delta=\frac{1}{2}$ and $t=1$, and multiplying by $A_n$.
Proposition \ref{MIR} may be checked numerically against the data
provided in Section \ref{sec:numerics}.

\subsection{Towards a Fredholm type determinant representation}
Expression \eqref{eq:efpMIR} is rather nice and compact. However, the
evaluation of its asymptotic behaviour for large $n$, $s$ is still out
of reach. For such a purpose, a Fredholm determinant
representation would be highly desiderable, in view of
the various well-established techniques developed in the last decades.

The notion of Fredhom determinant generalizes that of determinant of a
finite dimensional linear operator. It is defined for bounded
operators on a Hilbert space which differ from the identity operator
by a trace-class operator \cite{GGK-00}. In order to build a Fredholm
determinant representation, following \cite{TW-08b}, we have first
tried to rewrite the above integral representation as a Fredholm-type
determinant, that is as the determinant of some $s\times s$ matrix of
the form $1-M$.

We have managed to do so for the case $s=1,\dots, 4$.  This has
allowed us to guess a general structure for the matrix $M$, reported
in \eqref{eq:matrix}. The general structure is constrained by two
lemmas which can be easily proven, and by some natural recursive
properties, see \cite{CP-24} for details. The resulting formula,
Conjecture \ref{main_conj}, has gone through all our checks, see
below.

\section{Numerical checks}\label{sec:numerics}

\subsection{Monotone triangles}
To check Conjecture \ref{main_conj}, we want here to evaluate the
numbers $B_{n,s}$ by direct enumeration of ASMs with a frozen
square. We use the well-known bijection between ASMs and monotone
triangles \cite{MRR-83}.  A monotone triangle of size $n$ is a
triangular array of numbers $t_{ij}$, $i=1,\dots, n$, $j=1,\dots, i$,
such that:
  \begin{enumerate}
    \item for any given $i$, the sequence $t_{ij}$ is strictly
      increasing with $j$: $t_{ij}<t_{i,j+1}$, $j=1,\dots, i-1$;
    \item the sequences in two consecutive rows are weakly
      interlacing: $t_{i+1,j}\leq t_{ij}\leq t_{i+1,j+1}$;
    \item the $n$th (bottom) row is: $t_{nj}=j$, $j=1,\dots,n$.
      \end{enumerate}

Recall the correspondence between ASMs and monotone triangles: given
an $n\times n$ ASM with entries $a_{ij}$, compute a new $n\times n$
matrix with entries $b_{ij}=\sum_{l=1}^ia_{ij}$. It follows from the
definition of ASMs that the $i$th row of this matrix has exactly $i$
entries equal to $1$, all others being $0$.  Denoting the position of
the $j$th $1$'s in the $i$th row by $t_{ij}$, we indeed obtain a
monotone triangle, see Fig.~\ref{fig-dwbc_asm_mt}. The above described
mapping is easily seen to be bijective.

\subsection{Frozen-corner ASMs and  constrained monotone triangles}
Let us now consider an $n\times n$ ASM constrained to have an $s\times
s$ square of zeroes in, say, the top-right corner (recall that the
choice of the corner is irrelevant for the definition of $B_{n,s}$).
It is clear from the above bijection that such a condition translates
into the requirement that in the first $s$ rows of the corresponding
monotone triangle no entry is larger than $n-s$.  In other words,
$B_{n,s}$ is equal to the number of monotone triangle of size $n$,
satisfying the additional condition
\begin{equation}\label{eq:mt_frozen}
  t_{ij}\leq n-s,\qquad i\leq s.
  \end{equation}

We have implemented on Wolfram Mathematica \cite{W} the evaluation
of the number of monotone triangle satisfying condition
\eqref{eq:mt_frozen}. The code is a minor adaptation of that for
enumerating unconstrained monotone triangles, provided in
\cite{Br-99}. It is reported in Appendix \ref{app:code}.

\subsection{Results and checks}
We have run the code described above on a standard desktop and
evaluated $B_{n,s}$ for all $n\leq 20$, with $s\leq n$. We report
below the first few ones, till $n=12$.
\begin{longtable}[l]{rrl}
  $n$ & $s$ & $B_{n,s}$\rowspace
2 & 1 & \num{1}\rowspace
3 & 1 & \num{5}\rowspace
4 & 1 & \num{35}\\
4 & 2 & \num{4}\rowspace
5 & 1 & \num{387}\\
5 & 2 & \num{102}\rowspace
6 & 1 & \num{7007}\\
6 & 2 & \num{2889}\\
6 & 3 & \num{49}\rowspace
7 & 1 & \num{210912}\\
7 & 2 & \num{115089}\\
7 & 3 & \num{7007}\rowspace
8 & 1 & \num{10631868}\\
8 & 2 & \num{6994572}\\
8 & 3 & \num{901849}\\
8 & 4 & \num{1764}\rowspace
9 & 1 & \num{900985244}\\
9 & 2 & \num{673987304}\\
9 & 3 & \num{144298856}\\
9 & 4 & \num{1489302}\rowspace
10 & 1 & \num{128622437240}\\
10 & 2 & \num{105207867496}\\
10 & 3 & \num{32490142348}\\
10 & 4 & \num{945848476}\\
10 & 5 & \num{184041}\rowspace
11 & 1 & \num{30966210579675}\\
11 & 2 & \num{26959351496635}\\
11 & 3 & \num{10955022440189}\\
11 & 4 & \num{662654919034}\\
11 & 5 & \num{944518861}\rowspace
12 & 1 & \num{12580216114825125}\\
12 & 2 & \num{11439460153178700}\\
12 & 3 & \num{5734064378518415}\\
12 & 4 & \num{600268797504481}\\
12 & 5 & \num{3128505277443}\\
12 & 6 & \num{55294096}\\
\end{longtable}
\noindent
We have omitted vanishing values, occuring for $s > \lfloor
n/2 \rfloor $.  Further enumerations, for $n=13,\dots, 20$ are
reported in appendix \ref{app:numerics}.

 \begin{remark}
The numbers $B_{n,s}$ do not seem to factorize into `small' integers.
  \end{remark}
We have checked that the obtained values satisfy
properties \eqref{eq:property1}, \eqref{eq:property2},
\eqref{eq:property3}, \eqref{eq:property4} for all $n\leq 20$.  We
have also checked that the determinant formula proposed in Conjecture
\ref{main_conj} reproduces the obtained numbers for all $n\leq 20$ and
$s\leq n$.

 \section{Large ASMs, limit shape, and Tracy--Widom
 distribution}\label{sec:context}
 
\subsection{Arctic curve}
As said, the motivation in studying frozen-corner enumeration of ASMs
originates from asymptotic combinatorics, and in particular, from the
observed limit-shape phenomenon in large ASMs.  To be more precise,
for any given ASM $(a_{ij})_{1\leq i,j\leq n}$, we say that $(i,j)\in
[1,n]\times[1,n]$ is in the top-left frozen region if $a_{kl}=0$ for
all $k\leq i$ an $l\leq j$. We may similarly define the three other
frozen regions. If $(i,j)$ is not in any frozen region, then we say
that $(i,j)$ is in the disordered region.

Consider now the square lattice graph underlying any ASM; let
$\Lambda$ be its dual graph. The polygon living on $\Lambda$ and
separating the frozen and disordered regions is uniquely defined.  It
is called the \emph{frozen boundary}.  Clearly each given ASM will
have a given frozen boundary. A natural question concerns the
distribution of the frozen boundaries, when considering uniformly
sampled ASMs of given size $n$, and in particular if such distribution
concentrates, and where, in the limit $n\to\infty$.

Let us introduce the rescaled variables $x$ and $y$, with
$(x,y)\in[0,1]\times[0,1]$, related to the matrix indices $i,j$ as
follows, $i=\lceil y n\rceil$, $j=\lceil x n\rceil$.

\begin{theorem}\label{arctic}
  In the limit $n\to\infty$, the frozen boundary concentrate on a
  smooth curve (arctic curve) lying in the unit square.  Restricting
  to the top left quadrant, $x,y\in[0,1/2]$ for simplicity, the
  expression of the arctic curve in parametric form reads
\begin{equation}\label{eq:arctic}
  x=1-\frac{2\omega-1}{2\sqrt{\omega^2-\omega+1}}, \qquad
  y=1-\frac{\omega+1}{2\sqrt{\omega^2-\omega+1}}, \qquad \omega\in[1,\infty).
\end{equation}
The expression of the arctic curve in the other three quadrants may be
easily obtained from symmetry considerations.
\end{theorem}

\begin{remark}
The portion of arctic curve \eqref{eq:arctic} is the arc of ellipse of
equation:
\begin{equation}\label{eq:ellipse}
4x^2+4y^2-4xy-4x-4y+1=0,\qquad x,y\in\left[0,\textstyle\frac{1}{2}\right].
\end{equation}
\end{remark}

Note that the four points of tangency to the boundary of the unit
square are singular (discontinuous in the second
derivative). Expression \eqref{eq:arctic} has been derived by the
authors \cite{CP-08}, however one crucial step in the derivation
relied on some conjectural argument.  Subsequent developments, in
particular, the `tangent method' \cite{CS-16}, have opened the way for a
rigorous proof of Theorem \ref{arctic}, finally provided by
Aggarwal \cite{A-20}.

\subsection{Fluctuations around the arctic curve}
Let us turn to the fluctuations of the frozen boundary. A full
quantitative description of these is still beyond our
capabilities. Let us therefore restrict to the
\emph{intersection}, $(\xi,\xi)$, $\xi\in\{1,\dots,\lfloor n/2 \rfloor\}$,
of the frozen boundary with the ASM main diagonal. Under uniform
sampling of the ASMs of size $n$, the quantity $\xi$ is a discrete
random variable.  Its cumulative distribution function is simply
\begin{equation}
\mathbb{P}_n(\xi>s) = \frac{B_{n,s}}{A_n}.
\end{equation}
A relevant question is the behaviour of $\mathbb{P}_n(\xi> s)$ in
the limit of large $n$, $s$, with $s=\lfloor y n \rfloor$. However, to
address it, one needs a suitable expression for $B_{n,s}$, amenable to
asymptotic analysis. Remarkably, the formula proposed in Conjecture
\ref{main_conj} complies with such requirement.

Indeed, the Fredholm-type determinant may be equivalently rewritten as
the Fredholm determinant of some linear integral operator, with all
dependence on parameters $s$ and $n$ encoded in the kernel of the
integral operator.  The asymptotic behaviour of the kernel as a
function of the variable $y$ may then be evaluated by using the
saddle-point method.

In particular, in the vicinity of $\yc:=1-\sqrt{3}/2$, that is the
value of $y$ for which the arctic curve \eqref{eq:arctic} intersects
the main diagonal $x=y$, one has the coalescence of two
saddle-points. After suitable centering and rescaling of the variable
$y$ to $\sigma:=2^{4/3}3^{1/6}n^{2/3}(\yc-y)$, a careful analysis
shows that the kernel of the integral operator turns into the Airy
kernel, and thus, the corresponding Fredholm determinant into the
celebrated GUE Tracy--Widom distribution $\mathcal{F}_2(\sigma)$
\cite{TW-94a}. More precisely, the following result holds.
\begin{theorem}\label{TW}
  In the large $n$ limit, one has
\begin{equation}\label{eq:TW}
  \lim_{n\to\infty}\left[ \det_{1\leq i,j\leq
  s}(1-M)\Big|_{s=\big\lfloor
  n\yc-\frac{n^{1/3}}{2^{4/3}3^{1/6}}\sigma\big\rfloor}\right]
  =\mathcal{F}_2(\sigma),
\end{equation}
where $M$ is the $s\times s$ matrix entering the right hand side of
Conjecture \ref{main_conj}.
\end{theorem}
For full details on the derivation of the above result, see
\cite{CP-24}.
    
\subsection{Conclusion}
If Conjecture \ref{main_conj} is true, then Theorem \ref{TW}
implies that the fluctuations of the frozen boundary are governed by
Tracy--Widom distribution. Thus, taken together, these two statements
provide strong additional support to the occurrence of Tracy--Widom
distribution in large ASMs, already conjectured in \cite{ACJ-23}, and
checked against numerics in \cite{LKV-23,PS-23}.

Tracy--Widom distribution has already appeared in a variety of models
of relevance in combinatorics (see, e.g., \cite{D-06}, and reference
therein) or statistical mechanics (see, e.g., \cite{FS-24}, and
references therein). However, most of these models are `determinantal'
(or free-fermionic, in physics), while ASMs are not. Our result thus
extends the notion of `universality' of Tracy--Widom distribution, in
the sense of \cite{D-06}, to a wider class of models.

In conclusion, a rigorous derivation of Conjecture \ref{main_conj},
besides being of interest by itself, would fill a gap in the proof
that fluctuations of the frozen boundary of large ASMs are governed by
Tracy--Widom distribution, which has been a long-standing open
question.

\subsection*{Acknowledgements}
The authors are indebted to L. Cantini, S. Chhita, H. Spohn, and
J. Viti, for stimulating discussions. The work of A. G. Pronko was
supported by the Ministry of Science and Higher Education of the
Russian Federation (agreement 075-15-2025-344 dated 29/04/2025 for
Saint Petersburg Leonhard Euler International Mathematical Institute
at PDMI RAS).  The authors are grateful to an anonymous referee for
helpful remarks that improved the presentation.

\appendix

\section{}\label{app:B_proof}
We sketch here the derivation of identity \eqref{eq:property4}.  This
is much easier when resorting to the bijection between ASMs and path
configurations, shortly described in Section 1. 

\begin{figure}
\centering

\begin{tikzpicture}[scale=.5]
\draw [help lines] (0,1)--(10,1);
\draw [help lines] (0,2)--(10,2);
\draw [help lines] (0,3)--(10,3);
\draw [help lines] (0,4)--(10,4);
\draw [help lines] (0,5)--(10,5);
\draw [help lines] (0,6)--(10,6);
\draw [help lines] (0,7)--(10,7);
\draw [help lines] (0,8)--(10,8);
\draw [help lines] (0,9)--(10,9);
\draw [help lines] (1,0)--(1,10);
\draw [help lines] (2,0)--(2,10);
\draw [help lines] (3,0)--(3,10);
\draw [help lines] (4,0)--(4,10);
\draw [help lines] (5,0)--(5,10);
\draw [help lines] (6,0)--(6,10);
\draw [help lines] (7,0)--(7,10);
\draw [help lines] (8,0)--(8,10);
\draw [help lines] (9,0)--(9,10);
\draw [very thick](1,9)--(1,10);
\draw [very thick](2,9)--(2,10);
\draw [very thick](3,9)--(3,10);
\draw [very thick](4,9)--(4,10);
\draw [very thick](5,9)--(5,10);
\draw [very thick](6,9)--(6,10);
\draw [very thick](7,9)--(7,10);
\draw [very thick](8,9)--(8,10);
\draw [very thick](9,9)--(9,10);
\draw [very thick](0,1)--(1,1);
\draw [very thick](0,2)--(1,2);
\draw [very thick](0,3)--(1,3);
\draw [very thick](0,4)--(1,4);
\draw [very thick](0,5)--(1,5);
\draw [very thick](0,6)--(1,6);
\draw [very thick](0,7)--(1,7);
\draw [very thick](0,8)--(1,8);
\draw [very thick](0,9)--(1,9);
\draw [very thick] (1,9)--(5,9);
\draw [very thick] (1,8)--(5,8);
\draw [very thick] (1,7)--(5,7);
\draw [very thick] (1,6)--(5,6);
\draw [very thick] (1,5)--(1,9);
\draw [very thick] (2,5)--(2,9);
\draw [very thick] (3,5)--(3,9);
\draw [very thick] (4,5)--(4,9);
%
\draw [very thick] (7,5)--(7,6);
\draw [very thick] (4,2)--(5,2);
\fill [nearly transparent,gray] (4.2,.8)--(4.2,5.2)--(.8,5.2)--(.8,.8);
\fill [nearly transparent,gray] (4.8,5.8)--(9.2,5.8)--(9.2,9.2)--(4.8,9.2);
\fill [nearly transparent,gray] (4.8,.8)--(9.2,.8)--(9.2,5.2)--(4.8,5.2);
\draw [|->] (4.5,-.5)--(7.5,-.5);
\node at (6,-1) {$j$};
\draw [|->] (10.5,5.3)--(10.5,1.7);
\node at (11,3.5) {$k$};
\end{tikzpicture}

\caption{A $(2s+1)\times(2s+1)$ square ice grid, with an $s\times s$
frozen square in its top-left corner. The top-right and bottom-left
rectangles correspond to the refined enumerations of $s\times s$ ASMs
(see Fig.~\ref{fig-refined}), while the bottom-right $(s+1)\times
(s+1)$ square region contains only single path configurations.}
\label{fig-efp4}
\end{figure}

Let us consider a $(2s+1)\times(2s+1)$ ASM, with an $s\times s$ frozen
square in its top-left corner. In the osculating path picture, see
Fig.~\ref{fig-efp4}, this becomes a set of $(2s+1)$ path on the
$(2s+1)\times(2s+1)$ square grid, with an $s\times s$ square region in
the top-left corner where all vertices are of type 2.  Let us focus on
the $s\times (s+1)$ top-right portion of the grid: there are $s+1$
paths flowing from the top; $s$ of them have to exit from the left,
while the $(s+1)$th may only exit downward, say in correspondence of
the $(s+j)$th vertical edge, $j\in [1,\dots,s+1]$.  As explained above
see Fig.~\ref{fig-refined}, the number of such paths coincides with
the refined enumeration $A_{s+1,j}$ of ASMs.

Similarly, the number of path configurations in the $(s+1)\times s$
bottom-left portion of the grid, with one path entering horizontally
from the right at position $s+k$, coincides with the refined
enumeration $A_{s+1,k}$. Finally, for given $j$ and $k$, in the
$(s+1)\times (s+1)$ bottom-right portion of the grid we only have one
path, entering from the top, in correspondence of the $(s+j)$th
vertical edge, and exiting from the left, in correspondence of the
$(s+k)$th horizontal edge. The number of such single-path
configurations is simply $\binom{j+k-2}{j-1}$.  Multiplying these
three enumerations, and summing over the allowed values of $j$ and
$k$, immediately yields \eqref{eq:property4}.

\section{}\label{app:code}
We report here the Wolfram Mathematica code \cite{W} we have used
for the evaluation of the numbers $B_{n,s}$:
\begin{lstlisting}
nextlist[n_, alist_] :=
  Select[Distribute[Apply[
    Range, Partition[Prepend[Append[alist, n], 1], 2, 1], 1
      ], List], # == Union[#] &];
nextlistnew[n_, s_, alist_] := 
  If[Length[alist] < s, nextlist[n - s, alist],
    nextlist[n, alist]]
MT[n_, s_][blist_] := 
  MT[n, s][blist] = 
    If[Length[blist] >= n - 1, 1, 
      Apply[Plus, Map[MT[n, s], nextlistnew[n, s, blist]]
        ], {1}];
ASMfrozen[n_, s_] := Sum[MT[n, s][{i}], {i, n - s}]
\end{lstlisting}
This is just a minor adaptation of the code for unconstrained monotone
triangles provided in \cite{Br-99}, see page 63 therein.

\section{}\label{app:numerics}
We report here the number $B_{n,s}$, evaluated in terms of monotones
triangles, for $n=10,\dots,20$, and $s=1,\dots,\lfloor n /2\rfloor$.

\medskip

\begin{longtable}[l]{rrl}
  $n$ & $s$ & $B_{n,s}$\rowspace
13 & 1 & \num{8626772206437975000}\\
13 & 2 & \num{8085284537286414375}\\
13 & 3 & \num{4766467800019182375}\\
13 & 4 & \num{764765141679270475}\\
13 & 5 & \num{10238704485605765}\\
13 & 6 & \num{1757215525000}\rowspace
14 & 1 & \num{9986901971929869829500}\\  
14 & 2 & \num{9557633215532539485000}\\
14 & 3 & \num{6390208799423434022625}\\
14 & 4 & \num{1441330447360363147200}\\
14 & 5 & \num{39883284449338622740}\\
14 & 6 & \num{31549291544854616}\\
14 & 7 & \num{47675849104}\rowspace
15 & 1 & \num{19519080693305828937415200}\\
15 & 2 & \num{18949956228921855982385700}\\
15 & 3 & \num{13969410998712641496783000}\\
15 & 4 & \num{4152493908912199929945075}\\
15 & 5 & \num{204906100267058283679500}\\
15 & 6 & \num{506108207011474868400}\\
15 & 7 & \num{9503245338765360}\rowspace
16 & 1 & \num{64407656627191028079791846640}\\
16 & 2 & \num{63144359138386747743715394640}\\
16 & 3 & \num{50198563563919013122322373540}\\
16 & 4 & \num{18705823805304944943701823900}\\
16 & 5 & \num{1479675893764040036587431315}\\
16 & 6 & \num{8932759443524965167457800}\\
16 & 7 & \num{951434093994335362640}\\
16 & 8 & \num{117727187246656}\rowspace
17 & 1 & \num{358804774730434175758129259673456}\\
17 & 2 & \num{354105577359757906766371948697376}\\
17 & 3 & \num{298333903708787500221857600739936}\\
17 & 4 & \num{133912858837890294182713812604536}\\
17 & 5 & \num{15669580544590486163684302616040}\\
17 & 6 & \num{195234062777004744115907889066}\\
17 & 7 & \num{77617468143851451332322960}\\
17 & 8 & \num{148590087937175420256}\rowspace
18 & 1 & \num{3374501770056834424668067005366508160}\\
18 & 2 & \num{3345186415262559017533955220394878080}\\
18 & 3 & \num{2946094160268159399942635677670006000}\\
18 & 4 & \num{1542051040512127018766643769731903632}\\
18 & 5 & \num{250711193443677493960123784595089380}\\
18 & 6 & \num{5695769733133631831672720731794252}\\
18 & 7 & \num{6469250261084028486818521377683}\\
18 & 8 & \num{84778099221899163096886968}\\
18 & 9 & \num{831443906113411600}\rowspace
19 & 1 & \num{53576975973345090138316376635217569387788}\\
19 & 2 & \num{53270072499641082201198634936685040618732}\\
19 & 3 & \num{48517940677403477178301346505128327029748}\\
19 & 4 & \num{28829134880775563145866710298765418642088}\\
19 & 5 & \num{6193973979899248680612126583699721493780}\\
19 & 6 & \num{233438991061758469376565209482737854672}\\
19 & 7 & \num{628282787528898950170898927382130672}\\
19 & 8 & \num{36208603232194160709384879994300}\\
19 & 9 & \num{6693451140998292061585400}\rowspace
20 & 1 & \num{1435985354364704216564266803990901495284513968}\\
20 & 2 & \num{1430590451354449063721068755633229789198660380}\\
20 & 3 & \num{1336289815208055474404048759131669563982545640}\\
20 & 4 & \num{881410992082437335865516641683862633616426168}\\
20 & 5 & \num{240194284338153882283352726137805109829801860}\\
20 & 6 & \num{13940838952484992304062323074600042578016516}\\
20 & 7 & \num{77412505786986088748219280737825798922720}\\
20 & 8 & \num{14735602647010234326105847627175396848}\\
20 & 9 & \num{22142431679094238038503945163100}\\
20 & 10&  \num{16779127803917965290000}
\end{longtable}


\bibliography{nofile}

\end{document}